\documentclass[a4paper,english,12pt]{article}
\usepackage[logonly]{trace}
\usepackage{babel}
\usepackage{amsfonts, amsmath, amssymb}
\usepackage{graphicx}
\usepackage{pstricks}
\usepackage{psfig}
\usepackage{multirow}
\usepackage{fancyhdr}
\usepackage{vmargin, fancybox}

\newcommand{\mathsym}[1]{{}}

\begin{document}

\title{Spurious solitons and structural stability of finite difference schemes for nonlinear wave equations}

\author{Claire David\footnotemark[2] \footnote{Corresponding author:
david@lmm.jussieu.fr; fax number: (+33) 1.44.27.52.59.} ,  Pierre
Sagaut\footnotemark[2] \\ \small{Universit\'e Pierre et Marie Curie-Paris 6}  \\
\footnotemark[2] \small{ Laboratoire de Mod\'elisation en M\'ecanique, UMR CNRS 7607,} \\
\small{Bo\^ite courrier $n^0 162$, 4 place Jussieu, 75252 Paris,
cedex 05, France}}

%\date{December 2005}

\maketitle

\begin{abstract}
The goal of this work is to determine classes of traveling solitary
wave solutions for a differential approximation of a finite
difference scheme by means of a hyperbolic ansatz.
\end{abstract}

\section{Introduction}
\label{sec:intro}

\noindent The Burgers equation:

\begin{equation}
\label{Burgers}  u_t + c\, u\, u_x - \mu \,u_{xx} = 0,
\end{equation}

\noindent  $\alpha$, $\mu$ being real constants, plays a crucial
role in the history of wave equations. It was named after its use by
Burgers \cite{burger1} for studying
turbulence in 1939.\\

\noindent A finite difference scheme for the Burgers equation can be
written under the following form:

\begin{equation}
\label{scheme_general}  F({u_l}^m,h,\tau) = 0,
\end{equation}

\noindent where:
\begin{equation}
{u_l}^m=u\,(l\,dx, m\,dt)
\end{equation}
\noindent  $l\, \in \, \{i-1,\, i, \, i+1\}$, $m \, \in \, \{n-1,\,
n, \, n+1\}$, $j=1, \, ..., \, n_x$, $n=1, \, ..., \, n_t$, $h$,
$\tau$ denoting respectively the mesh size and time step, and
$\sigma$
the Courant-Friedrichs-Lewy number ($cfl$) coefficient, defined as $\sigma = c \,\tau / h$.\\

\noindent A numerical scheme is  specified by selecting appropriate
expression of the function $F$ in equation (\ref{scheme_general}).\\

\noindent Considering the ${u_l}^m$ terms as functions of the mesh
size $h$ and time step $\tau$, expanding them at a given order by
means of their Taylor series, and neglecting the $o(\tau^p)$ and
$o({h}^q)$ terms, for given values of the integers $p$, $q$, lead to
a differential approximation of the Burgers equation, of the form:

\begin{equation}
\label{Approx}  \mathcal {F} (u,\frac {\partial^r u}{\partial
x^r},\frac {\partial^s u}{\partial t^s},h,\tau) = 0,
\end{equation}
\noindent $r$, $s$ being integers.\\

\noindent For sake of simplicity, a non-dimensional form of
(\ref{Approx}) will be used:

\begin{equation}
\label{ApproxAdim}  \widetilde{\mathcal {F}} (\tilde{u},\frac
{\partial^r \tilde{u}}{\partial \tilde{x}^r},\frac {\partial^s
\tilde{u}}{\partial \tilde{t}^s}) = 0,
\end{equation}

\noindent Depending on this differential approximation
(\ref{Approx}), solutions,
as solitary waves, may arise.\\

\noindent The paper is organized as follows. Two specific schemes
are exhibited in section \ref{Section:Scheme}. The general method is
exposed in Section \ref{Th}. In Section \ref{Section:Sol}, it is
shown that out of the two studied schemes, only one leads to
solitary waves. A related class of traveling wave solutions of
equation (\ref{Approx}) is thus presented, by using a hyperbolic
ansatz. The stability of this class of solutions is discussed in
section \ref{Section:Stab}. A numerical example is exposed in
section \ref{Ex}.

\section{Scheme study}
\label{Section:Scheme}

\subsection{Finite second-order centered scheme in space,\\ Euler-time scheme} \label{Scheme1}

 \noindent For the finite second-order centered scheme in space and Euler-time scheme, the function $F$ of
 (\ref{scheme_general}) takes the form:

\begin{equation}
\label{ApproxEuler}  F({u_l}^m,h,\tau)=\frac
{{u_i}^{n+1}-{u_i}^{n}}{\tau}+ c\, {u_i}^{n}\,\frac
{{u_{i+1}}^{n}-{u_{i-1}}^{n}}{2\,h} -\mu \,\frac
{{u_{i+1}}^{n}-2\,{u_i}^{n}+{u_{i-1}}^{n}}{{h}^2}=0 \end{equation}

\noindent Consider ${u_{i}}^{n+1}$ as a function of the time step
$\tau$, and expand it at the second order by means of its Taylor
series:

\begin{equation}
{u_i}^{n+1}=u\,(i\,h, (n+1)\,\tau)=u\,(i\,h,
n\,\tau)+\tau\,u_t\,(i\,h,
n\,\tau)+\frac{{\tau}^2}{2}\,u_{tt}\,(i\,h, n\,\tau)+o({\tau}^2)
\end{equation}

\noindent It ensures:

\begin{equation}
\label{Prop1} \frac{{u_i}^{n+1}-{u_i}^{n}}{\tau}=u_t\,(i\,h,
n\,\tau)+\frac{{\tau}}{2}\,u_{tt}\,(i\,h, n\,\tau)+o({\tau})
\end{equation}

\noindent In the same way, consider ${u_{i+1}}^{n}$ and
${u_{i-1}}^{n}$ as functions of the mesh size $h$, and expand them
at the fourth order by means of their Taylor series:
\begin{equation}
\label{Prop_a} \scriptsize{ \begin{aligned}
{u_{i+1}}^{n}  =   & u\,((i+1)\,h, n\,\tau)\\
=   &  u\,(i\,h, n\,\tau)+h\,u_x\,(i\,h,
n\,\tau)+\frac{{h}^2}{2}\,u_{xx}\,(i\,h,
n\,\tau)+\frac{{h}^3}{3!}\,u_{xxx}\,(i\,h,
n\,\tau)+\frac{{h}^4}{4!}\,u_{xxxx}\,(i\,h, n\,\tau)+o({h}^4)
\end{aligned}}
\end{equation}

\begin{equation}
\label{Prop_b} \scriptsize{ \begin{aligned}
{u_{i-1}}^{n}  =   & u\,((i-1)\,h, n\,\tau)\\
=   &  u\,(i\,h, n\,\tau)-h\,u_x\,(i\,h,
n\,\tau)+\frac{{h}^2}{2}\,u_{xx}\,(i\,h,
n\,\tau)-\frac{{h}^3}{3!}\,u_{xxx}\,(i\,h,
n\,\tau)+\frac{{h}^4}{4!}\,u_{xxxx}\,(i\,h, n\,\tau)+o({h}^4)
\end{aligned}}
\end{equation}

\noindent It ensures:

\begin{equation}
\label{Prop2} \frac {{u_{i+1}}^{n}-2\,{u_i}^{n}+ {u_{i-1}}^{n}}
{{h}^2}=u_{xx}\,(i\,h,
n\,\tau)+\frac{{2\,h}^2}{4!}\,u_{xxxx}\,(i\,h, n\,\tau)+o({h}^2)
\end{equation}

\noindent and:

\begin{equation}
\label{Prop3} \frac{{u_{i+1}}^{n}- {u_{i-1}}^{n}}{2\,h}=u_x\,(i\,h,
n\,\tau)+\frac{{h}^2}{3!}\,u_{xxx}\,(i\,h, n\,\tau)+o({h}^3)
\end{equation}

\bigskip \noindent  Equation (\ref{Approx}) can thus be written as:

\begin{equation}
\label{BurgersEq} \scriptsize{ \begin{aligned} & u_t\,(i\,h,
n\,\tau)+\frac{{\tau}}{2}\,u_{tt}\,(i\,h, n\,\tau)+o({\tau})\\ &+
 c\, u\,(i\,h, n\,\tau) \, \big [ u_x\,(i\,h,
n\,\tau)+\frac{{h}^2}{3!}\,u_{xxx}\,(i\,h, n\,\tau)+o({h}^3)\,\big ]\\
&-\mu \,\big [ u_{xx}\,(i\,h,
n\,\tau)+\frac{{2\,h}^2}{4}\,u_{xxxx}\,(i\,h, n\,\tau)+o({h}^2) \big
] =0
\end{aligned}}
\end{equation}

\noindent i. e., at $x=i\,h$ and $t=n\,\tau$:

\begin{equation}
 \scriptsize{  \big [\,
u_t+\frac{{\tau}}{2}\,u_{tt}+o({\tau})+  c\, u\,  \big [
u_x+\frac{{h}^2}{3!}\,u_{xxx}+o({h}^3) \big ]-\mu \,\big [
u_{xx}+\frac{{2\,h}^2}{4!}\,u_{xxxx}+o({h}^2) \big ] \, \big
]_{(x,t)}=0}
\end{equation}

\noindent The first differential approximation of the Burgers
equation (\ref{Burgers}) is thus obtained neglecting the $o(\tau)$
and $o({h}^2)$ terms:

\begin{equation}
 \scriptsize{  \big [\,
u_t+\frac{{\tau}}{2}\,u_{tt}+ c\, u\, \big [
u_x+\frac{{h}^2}{3!}\,u_{xxx}  \big ]-\mu \,\big [
u_{xx}+\frac{{h}^2}{12}\,u_{xxxx} \big ]\, \big ]_{(x,t)}=0}
\end{equation}

\noindent that we will keep as:

\begin{equation}
\label{Eq} u_t+ c\, u\,u_x-\mu \,
u_{xx}+\frac{{\tau}}{2}\,u_{tt}+\frac{{h}^2}{6}\,u\,u_{xxx}-\mu
\,\frac{{h}^2}{12}\,u_{xxxx}  =0
\end{equation}

\noindent For sake of simplicity, this latter equation can be
adimensionalized through in the following way:\\
\noindent set:

\begin{equation}
\label{Chgt_Var_Adim} \left \lbrace
\begin{array}{rcl}
u &  =   & U_0\,\tilde {u}\\
t &  =   & \tau_0\,\tilde {t}\\
x &  =   & h_0\,\tilde {x}
\end{array}
\right. \end{equation}

\noindent where:

\begin{equation}
\label{U0} U_0=\frac {h_0}{\tau_0}\end{equation}

\noindent Denote by $ Re_h$ the mesh Reynolds number, defined as:

\begin{equation}
\label{Re}
 Re_h=\frac{U_0\,h}{\mu}
\end{equation}

\noindent The change of variables (\ref{Chgt_Var_Adim}) leads to:

\begin{equation}
\label{Chgt_Var} \left \lbrace
\begin{array}{rcl}
u_t &  =   & \frac {U_0}{\tau_0} \,\tilde {u}_{\tilde {t}}\\
u_{x^k} &  =   & \frac {U_0}{h_0^k} \,\tilde {u}_{{\tilde {x}}^k }
\end{array}
\right. \end{equation}

\noindent Multiplying (\ref{Eq}) by $\frac {\tau_0}{U_0}$ yields:

\begin{equation}
 \tilde {u}_{\tilde {t}}+ c\, \frac{U_0\,\tau_0}{h_0}\, \tilde {u}\, \tilde {u}_{\tilde {x}}-\mu \,
\, \frac{\tau_0}{h_0^2}\, {\tilde {u}}_{\tilde {x}\tilde
{x}}+\frac{{\tau}}{2\,\tau_0}\,\tilde {u}_{\tilde {t}\tilde
{t}}+\frac{{h}^2\,U_0\tau_0}{6\,h_0^3}\,\tilde {u}\,\tilde
{u}_{\tilde {x}\tilde {x}\tilde {x}}-\mu
\,\frac{{h}^2\,\tau_0}{12\,h_0^4}\,\tilde {u}_{\tilde {x}\tilde
{x}\tilde {x}\tilde {x}} =0
\end{equation}

\noindent Relations (\ref{U0}) and (\ref{Re}) ensure:

\begin{equation}
\label{inter}
 \tilde {u}_{\tilde {t}}+ c\, \tilde {u}\, \tilde {u}_{\tilde {x}}-\frac{h}{h_0\,Re_h}
 \, {\tilde {u}}_{\tilde {x}\tilde
{x}}+\frac{{\tau}}{2\,\tau_0}\,\tilde {u}_{\tilde {t}\tilde
{t}}+\frac{{h}^2}{6\,h_0^2}\,\tilde {u}\,\tilde {u}_{\tilde
{x}\tilde {x}\tilde {x}}-\frac{{h}^3}{12\,Re_h\,h_0^3}\,\tilde
{u}_{\tilde {x}\tilde {x}\tilde {x}\tilde {x}} =0
\end{equation}

\noindent For $h=h_0$, due to $\sigma=\frac {U_0\,\tau}{h}$,
(\ref{inter}) becomes:

\begin{equation}
\label{inter}
 \tilde {u}_{\tilde {t}}+ c\, \tilde {u}\, \tilde {u}_{\tilde {x}}-\frac{1}{Re_h}
 \, {\tilde {u}}_{\tilde {x}\tilde
{x}}+\sigma\,\tilde {u}_{\tilde {t}\tilde {t}}+\frac{1}{6}\,\tilde
{u}\,\tilde {u}_{\tilde {x}\tilde {x}\tilde
{x}}-\frac{1}{12\,Re_h}\,\tilde {u}_{\tilde {x}\tilde {x}\tilde
{x}\tilde {x}} =0
\end{equation}

\subsection{The Lax-Wendroff scheme}

\noindent For the Lax-Wendroff scheme, the function $F$ of
 (\ref{scheme_general}) takes the form:

\begin{equation}
\label{ApproxLW}  F({u_l}^m,h,\tau)=\frac
{{u_i}^{n+1}-{u_i}^{n}}{\tau}+ c \,{u_i}^{n}\, \left \lbrace \frac
{{u_{i+1}}^{n}-{u_{i-1}}^{n}}{2\,h}\right \rbrace -(\mu +\frac {c^2
\,\tau} {2})\,\left \lbrace \frac
{{u_{i+1}}^{n}-2\,{u_i}^{n}+{u_{i-1}}^{n}}{{h}^2}\right \rbrace =0
\end{equation}

\noindent $\frac{{u_i}^{n+1}-{u_i}^{n}}{\tau}$ is expressed by means
of (\ref{Prop1}), and $\frac {{u_{i+1}}^{n}-2\,{u_i}^{n}+
{u_{i-1}}^{n}} {{h}^2}$ by means of (\ref{Prop2}):

\begin{equation}
\frac {{u_{i+1}}^{n}-2\,{u_i}^{n}+ {u_{i-1}}^{n}}
{{h}^2}=u_{xx}\,(i\,h,
n\,\tau)+\frac{{2\,h}^2}{4!}\,u_{xxxx}\,(i\,h, n\,\tau)+o({h}^2)
\end{equation}

\noindent (\ref{Prop2}) also yields:

\begin{equation}
\frac{{u_{i+1}}^{n}- {u_{i-1}}^{n}}{2\,h}=u_x\,(i\,h,
n\,\tau)+\frac{{h}^2}{3!}\,u_{xxx}\,(i\,h, n\,\tau)+o({h}^3)
\end{equation}

\bigskip \noindent Equation (\ref{ApproxLW}) can thus be written as:

\begin{equation}
\label{BurgersEq} \scriptsize{ \begin{aligned} & u_t\,(i\,h,
n\,\tau)+\frac{{\tau}}{2}\,u_{tt}\,(i\,h, n\,\tau)+o({\tau})\\ &+
 \alpha\, u\,(i\,h, n\,\tau) \,\big [ u_x\,(i\,h,
n\,\tau)+\frac{{h}^2}{3!}\,u_{xxx}\,(i\,h, n\,\tau)+o({h}^3)\,\big ]\\
&-(\mu +\frac {c^2 \,\tau} {2}) \,\big [ u_{xx}\,(i\,h,
n\,\tau)+\frac{{2\,h}^2}{4!}\,u_{xxxx}\,(i\,h, n\,\tau)+o({h}^2)
\big ] =0
\end{aligned}}
\end{equation}

\noindent i. e., at $x=i\,h$ and $t=n\,\tau$:

\begin{equation}
 \scriptsize{  \big [\,
u_t+\frac{{\tau}}{2}\,u_{tt}+o({\tau})+  c\, u\, \big [
u_x+\frac{{h}^2}{3!}\,u_{xxx}+o({h}^3) \big ]-(\mu +\frac {c^2
\,\tau} {2}) \,\big [ u_{xx}+\frac{{2\,h}^2}{4!}\,u_{xxxx}+o({h}^2)
\big ] \, \big ]_{(x,t)}=0}
\end{equation}

\noindent The first differential approximation of the Burgers
equation (\ref{Burgers}) is thus obtained neglecting the $o(\tau)$
and $o({h}^2)$ terms:

\begin{equation}
 \scriptsize{  \big [\,
u_t+\frac{{\tau}}{2}\,u_{tt}+ c\, u\,\big [\,
u_x+\frac{{h}^2}{3!}\,u_{xxx}\, \big ]-(\mu +\frac {c^2 \,\tau} {2})
\,\big [ u_{xx}+\frac{{h}^2}{12}\,u_{xxxx} \big ] \, \big
]_{(x,t)}=0}
\end{equation}

\noindent that we will keep as:

\begin{equation}
\label{Eq1LW} u_t+ c\, u\,u_x-(\mu +\frac {c^2 \,\tau} {2\,h^2}) \,
u_{xx}+\frac{{\tau}}{2}\,u_{tt}+\frac{{h}^2}{6}\,u\,u_{xxx}-(\mu
+\frac {c^2 \,\tau} {2}) \,\frac{{h}^2}{12}\,u_{xxxx}  =0
\end{equation}

\noindent Equation (\ref{Eq1LW}) is adimensionalized as in section
\ref{Scheme1}, leading to:

\begin{equation}
\label{Eq1LW_adim} \tilde{u}_{\tilde{t}}+ c\,
\tilde{u}\,\tilde{u}_x-(\frac{1}{Re_h} +\frac {c^2 \,\sigma} {2}) \,
\tilde{u}_{\tilde{x}\tilde{x}}+\frac{{\tau}}{2}\,\tilde{u}_{\tilde{t}\tilde{t}}
+\frac{1}{6}\,\tilde{u}\,\tilde{u}_{\tilde{x}\tilde{x}\tilde{x}}-(\frac{1}{Re_h}
+\frac {c^2 \,\sigma} {2})
\,\frac{1}{12}\,\tilde{u}_{\tilde{x}\tilde{x}\tilde{x}\tilde{x}} =0
\end{equation}

\section{Solitary waves}
\label{Th}

\noindent Approximated solutions of the Burgers equation
(\ref{Burgers}) by means of the difference scheme
(\ref{scheme_general}) strongly depend on the values of the time and
space steps. For specific values of $\tau$ and $h$, equation
(\ref{ApproxAdim}) can, for instance, have traveling wave solutions
which can be of great
disturbance to the searched solution.\\
\noindent We presently aim at determining the conditions, depending
on the values of the parameters $\tau$ and $h$, which give birth to
traveling wave solutions of (\ref{Eq}).\\
\noindent Following Feng \cite{feng1} and our previous work
\cite{David}, where traveling wave solutions of the cBKDV equation
were exhibited as combinations of bell-profile waves and
kink-profile waves, we aim
at determining traveling wave solutions of (\ref{ApproxAdim}).\\
\noindent Following \cite{feng1}, we assume that equation
(\ref{ApproxAdim}) has the traveling wave solution of the form
\begin{equation}\label{ChgtVar} \tilde{u}(\tilde{x}, \tilde{t}) =\tilde{u}(\xi), \quad \xi= \tilde{x}-v\,\tilde{t} \end{equation}
where $v$ is the wave velocity. Substituting (\ref{ChgtVar}) into
equation (\ref{ApproxAdim}) leads to:

\begin{equation}
\label{EqXi}  \widetilde{\mathcal {F}}_\xi
(\tilde{u},\tilde{u}^{(r)},(-v)^s\,\tilde{u}^{(s)}) = 0,
\end{equation}

\noindent Performing an integration of (\ref{EqXi}) with respect to
$\xi$ and setting the integration constant to zero leads to an
equation of the form:

\begin{equation}
\label{EqXiInt}  \widetilde{\mathcal {F}}_\xi^{\mathcal P}
(\tilde{u},\tilde{u}^{(r)},(-v)^s\,\tilde{u}^{(s)}) = 0,
\end{equation}

\noindent which will be the starting point for the determination of
solitary waves solutions.

\section{Traveling Solitary Waves}

\label{Section:Sol}

\subsection{Hyperbolic Ansatz}

\noindent The discussion in the preceding section provides us useful
information when we construct traveling solitary wave solutions for
equation (\ref{EqXi}). Based on this result, in this section, a
class of traveling wave solutions of the equivalent equation
(\ref{Eq}) are searched as a combination of bell-profile waves and
kink-profile waves of the form
\begin{equation}
\label{sol} \tilde{u}(\tilde{x}, \tilde{t}) = \sum_{i = 1}^n \left
(U_i\; \text{tanh}^i \left [\,C_i (\tilde{x}-v\,\tilde{t})\, \right
] + V_i \; \text{sech}^i \left [\,C_i(\tilde{x}-v\,\tilde{t}+x_0)\,
\right ] \right )+V_0
\end{equation}
where the $U_i's$, $V_i's$,  $C_i's$, $(i=1,\ \cdots,\ n)$, $V_0$
and $v$ are constants to be determined.
\\
\noindent In the following, $c$ is taken equal to 1.

\subsection{Theoretical study}

\noindent Substitution of (\ref{sol}) into equation (\ref{EqXi})
leads to an equation of the form

\begin{equation}
\label{Eqgen}
 \sum_{i,\, j, \,k} A_i\,\text{tanh}^i \big (C_i\,
\xi \big )\,\text{sech}^j \big (C_i \,\xi \big )\,\text{sinh}^k \big
(C_i\, \xi  \big )  =0
\end{equation}

\noindent the $A_i$ being real constants.\\
 \noindent The difficulty
for solving equation (\ref{Eqgen}) lies in finding the values of the
constants $U_i$, $V_i$,  $C_i$, $V_0$ and $v$ by using the
over-determined algebraic equations. Following \cite{feng1}, after
balancing the higher-order derivative term and the leading nonlinear
term, we deduce $n=1$. Then, following \cite{David} we replace
$\mbox{sech}({C_1} \,\xi)$ by $ \frac{2}{ e^{\,{C_1} \,\xi}+e^{\,-
{C_1} \,\xi } }$, $\mbox{sinh(}{C_1} \,\xi)$ by $ \frac { e^{\,{C_1}
\,\xi}-e^{\,- {C_1} \,\xi } }{2}$, $\mbox{tanh}({C_1} \,\xi)$ by
$\frac{e^{\, {C_1} \,\xi } -
      e^{\,- {C_1} \,\xi}}{e^{ \,{C_1} \,,\xi   } + e^{\,-{C_1} \,\xi }}$,  and multiply both sides by
      ${({  e^{\xi\, {C_1}}+e^{\,- \xi\,
{C_1}  } })}^{5}\, e^{\,5\,\xi \, {C_1}}$, so that equation
(\ref{Eqgen}) can be rewritten in the following form:
\begin {equation}
\label{Syst} \sum_{k=0}^{10} P_k ( U_1,\ V_1,\ C_1,\ v,\, V_0 )
\,e^{\,k \,C_1\, \xi} \; = \;0,
\end{equation}
where the $P_k$ $(k=0,\,...,\,10)$, are polynomials of $U_1$, $V_1$,
$C_1$, $V_0$  and $v$.

\subsection{Scheme study}

\subsubsection{The Finite second-order centered scheme in space, Euler-time scheme} \label{Scheme1_Sol}

\noindent Equation (\ref{EqXi}) is presently given by:

\begin{equation}
\label{EqXi1} -v\,\tilde{u}'( \xi )+ c\, \tilde{u}( \xi
)\,\tilde{u}'( \xi )-\frac{1}{Re_h} \, \tilde{u}''( \xi
)+v^2\,\frac{{\tau}}{2}\,\tilde{u}''( \xi )+\frac{1}{6}\,\tilde{u}(
\xi )\,\tilde{u}^{(3)}( \xi )-\frac{1}{Re_h}
\,\frac{1}{12}\,\tilde{u}^{(4)}( \xi )  =0
\end{equation}
\noindent Performing an integration of (\ref{EqXi1}) with respect to
$\xi$ and setting the integration constant to zero yields:

\begin{equation}
\label{EqXi2} -v\,\tilde{u}( \xi )+\frac{ c}{2}\,\tilde{u}^2( \xi
)+(v^2\,\frac{{\sigma}}{2}-\frac{1}{Re_h})\,{\tilde{u}}'( \xi )
+\frac{1}{6}\,\left \lbrace \tilde{u}(\xi)\,{\tilde{u}}^{(2)}( \xi
)-\frac{1}{2}\,{\tilde{u}}^{'^2}( \xi )\right
\rbrace-\frac{1}{12\,Re_h}\,\tilde{u}^{(3)}( \xi ) =0
\end{equation}

\noindent The related system (\ref{Syst}) has consistent solutions,
which are given in Tables \ref{Table1}.\\
\noindent For sake of simplicity, we use $\varepsilon$ to denote $1$
or $-1$.

\begin{table}[h!]
\caption{}
\begin{center}
\begin{tabular}{cccccccccc}
\hline
 & $\sigma$ & $v$ & $U_1 $  & $V_1$  & $C_1$& $V_0$ \\
\hline \hline Sets 1, 2 &  $\frac{484\,Re_h}{729}$ &
$\scriptsize{\varepsilon  \,\frac{108}{11 \sqrt{11}\,Re_h}}$ &
$\scriptsize{\varepsilon\,\frac{108}{5 \sqrt{11}\,Re_h}} $  & $0$  &
$\scriptsize{-\varepsilon\,\frac{6}{\sqrt{11}}}$&
$\scriptsize{-\varepsilon\,\frac{108}{5 \sqrt{11}\,Re_h}}$ \\
Set 3 & $\scriptsize{\frac{5 \,Re_h \,\left(17\,C_1^2-12
  \right)}{6 C_1^2 \left(4 C_1^2-9\right)^2}}$ & $\scriptsize{-\frac{2 \left(4 C_1^3-9 C_1\right)}{5\,Re_h}}$ &
   $\scriptsize{-\frac{18 C_1}{5\,Re_h} }$  & $0$  & $\scriptsize{ \in   {\mathbb R} } $& $\frac{18 C_1}{5\,Re_h}$ \\
Set 4   & $-\frac{\,Re_h \left(64 C_1^6-384 C_1^4+551
   C_1^2-156\right)}{6 C_1^2 \left(4 C_1^2-9\right)^2}$ & $\scriptsize{\frac{-\frac{5 C_1}{13-8 C_1^2}-C_1}{Re_h}}$ &
   $\scriptsize{-\frac{2 \left(8 C_1^3-9 C_1\right)}{\,Re_h \left(8
   C_1^2-13\right)} }$  & $0$  & $\scriptsize{ \in   {\mathbb R} } $ &
   $\scriptsize{\frac{18 C_1}{\,Re_h \left(8 C_1^2-13\right)}}$ \\
\hline
\end{tabular}
\end{center}
\label{Table1}
\end{table}%

\noindent In the following, we shall denote:

\begin{equation}
\label{Chgt_Var} \left \lbrace
\begin{array}{rcl}
\sigma_{1,2} &  = & \frac{484\,Re_h}{729}\\
\sigma_3 &  =   & \scriptsize{\frac{5 \,Re_h \,\left(17\,C_1^2-12
  \right)}{6 C_1^2 \left(4 C_1^2-9\right)^2}}\\
\sigma_4  &  =   & -\frac{\,Re_h \left(64 C_1^6-384 C_1^4+551
   C_1^2-156\right)}{6 C_1^2 \left(4 C_1^2-9\right)^2}\,=   \, -\frac{\,Re_h \,(C_1^2-4)\,(8\,C_1^2-13)\,(8\,C_1^2-3)\,}{6 C_1^2 \left(4 C_1^2-9\right)^2}
\end{array}
\right. \end{equation}

\subsubsection{The Lax-Wendroff scheme}

\noindent Equation (\ref{EqXi}) is then given by:

\begin{equation}
\label{EqXiLW} \scriptsize{-v\, \tilde {u}'( \xi )+ c\,  \tilde {u}(
\xi )\,u'( \xi )-\big(\frac{1}{Re_h}+\frac {c^2 \,\sigma} {2}\big)
\, \tilde {u}''( \xi )+v^2\,\frac{{\sigma}}{2}\, \tilde {u}''( \xi
)+\frac{1}{6}\, \tilde {u}( \xi )\,u^{(3)}( \xi )-(\frac{1}{Re_h}
+\frac {c^2 \,\sigma} {2}) \,\frac{1}{12}\,u^{(4)}( \xi )  =0}
\end{equation}

\noindent Performing an integration with respect to $\xi$ and
setting the integration constant to zero yields:

\begin{equation}
\label{EqXiLW2} \scriptsize{-v\, \tilde {u}( \xi )+\frac{ c}{2}\,
\tilde {u}^2( \xi )+(v^2\,\frac{{\sigma}}{2}-(\frac{1}{Re_h} +\frac
{c^2 \,\sigma} {2}))\,u'( \xi ) +\frac{1}{6}\,\big [\tilde
{u}(\xi)\, \tilde {u}^{(2)}( \xi )-\frac{1}{2}\,{  \tilde
{u}}^{'^2}( \xi )\big ]-\big(\frac{1}{Re_h} +\frac {c^2 \,\sigma}
{2}\big) \,\frac{1}{12}\, \tilde {u}^{(3)}( \xi ) =0}
\end{equation}

\noindent The related system (\ref{Syst}) does not admit consistent
solutions.

\section{Stability study}
 \label{Section:Stab}

\noindent In the following, the stability of the solutions presented
in section \ref{Scheme1_Sol} is discussed.\\

\noindent The variations of the $cfl$ coefficient $\sigma$ as a
function of the parameters $Re_h$, $C_1$, have crucial influence on
the stability.\\
\noindent $\sigma_3$, $\sigma_4$ being even functions of $C_1$, we
shall restrain our study to $C_1 \,\in \,]\,0,+\infty\,[$.\\
\noindent Calculation yield:

\begin{equation}
\left \lbrace
\begin{array}{rcl}
\frac{\partial {\frac {\sigma_3}{Re_h}}} {\partial
{C_1^2}}&=&-\frac{20 \left(34 \,{C_1}^4-36
   \,{C_1}^2+27\right) }{3 \,C_1^3
   \left(4 \,C_1^2-9\right)^3} \,> \,0\,\,\, ,\,\,\,\frac{\partial {\big [\sigma_3\cdot Re_h\big ]}} {\partial {C_1^2}}
   =-\frac{20 \left(34
\,{C_1}^4-36
   \,{C_1}^2+27\right) \,Re_h^2}{3 \,C_1^3
   \left(4 \,C_1^2-9\right)^3} \,> \,0\\
   \frac{\partial {\big [\frac {\sigma_4}{Re_h}\big ]}} {\partial {C_1^2}}&=&-\frac{4 \left(96 \,C_1^6-238 \,C_1^4+468
   \,C_1^2-351\right) \,Re_h}{3 \,C_1^3
   \left(4 \,C_1^2-9\right)^3}\,\,\, ,\,\,\,
   \frac{\partial {\big [\sigma_4\cdot Re_h\big ]}} {\partial {C_1^2}}=-\frac{4 \left(96 \,C_1^6-238 \,C_1^4+468
   \,C_1^2-351\right) \,Re_h^2}{3 \,C_1^3
   \left(4 \,C_1^2-9\right)^3}
\end{array}
\right. \end{equation}

\noindent Denote by ${C_1^0}^2$ the value of ${C_1}^2$ for which
$\frac{\partial {\big [\frac {\sigma_4}{Re_h}\big ]}} {\partial
{C_1}}$ and $ \frac{\partial {\big [\sigma_4\,Re_h\big ]}} {\partial
{C_1}}$ vanish.\\

\noindent For a given value of the mesh Reynolds number $Re_h$, we
obtain the following interesting variation tables:

\begin{equation}
\begin{array}{|c|ccccr|}
\hline
C_1^2    & 0  &    &  \frac{3}{2}  &     & +\infty  \\
\hline \frac{\partial {\big [\frac {\sigma_3}{Re_h}\big ]}}
{\partial {C_1^2}}\,,\,
\frac{\partial {\big [\sigma_3\cdot Re_h\big ]}} {\partial {C_1^2}} &          &  + & \|  & + &   \\
\hline
 & & &  +\infty \, \|\,  +\infty& & \\       % ligne des valeurs "max"
\frac {\sigma_3}{Re_h}\,,\, {\sigma_3}\cdot{Re_h}&  &\nearrow   &  &\searrow  & \\   % flèches
&  - \infty & & \|\,  & & 0 \\          % ligne des valeurs "min"
\hline
\end{array}
\end{equation}

\begin{equation}
\begin{array}{|c|ccccccr|}
\hline
C_1^2    & 0  &  & {{{C_1}^0}^2}  &  & \frac{3}{2}   &     & +\infty  \\
\hline \frac{\partial {\big [\frac {\sigma_4}{Re_h} \big ]}}
{\partial {C_1^2}}\,,\,\frac
{\partial {\big [\sigma_4\cdot Re_h\big ]}}{\partial {C_1^2}}  &      +    & &  0 & -  & \|  & + &   \\
\hline
 &+\infty & &   &&+ \infty\, \|\, + \infty  &  & \\       % ligne des valeurs "max"
\frac {\sigma_4}{Re_h} \,,\, {\sigma_4}\cdot{Re_h}&  &\searrow    & & \nearrow &&\searrow  & \\   % flèches
&  & &<0  && & & <0  \\          % ligne des valeurs "min"
\hline
\end{array}
\end{equation}

\noindent $\frac {\sigma}{Re_h}$, ${\sigma}\,{Re_h}$ take all the
values between 0 and $+\infty$. Thus, for stable and unstable sets
$\big ( {\sigma}\cdot {Re_h} \big )$, there exists a solitary wave
solution of (\ref{EqXi2}).

\section{Numerical Example}
\label{Ex}

\noindent In the following, we specifically consider the third
traveling solitary wave (see (\ref{Table1})) solution of equation
(\ref{EqXi2}).\\
\noindent Numerical values of the parameters are: $Re_h=1.9$,
$C_1=3$.

\subsection{Analytical soliton}

\noindent The third traveling solitary wave solitary wave (see
(\ref{Table1})) solution of equation (\ref{EqXi2}), is given by:

\begin{equation}
\label{Soliton} \tilde{u}(\tilde{x}, \tilde{t}) = U_1\; \text{tanh}
\left [\,C_1 \,(\tilde{x}-v\,\tilde{t})\, \right ] +V_0
\end{equation}

\begin{figure}[h!]
\hskip 4cm
 \includegraphics[width=7cm]{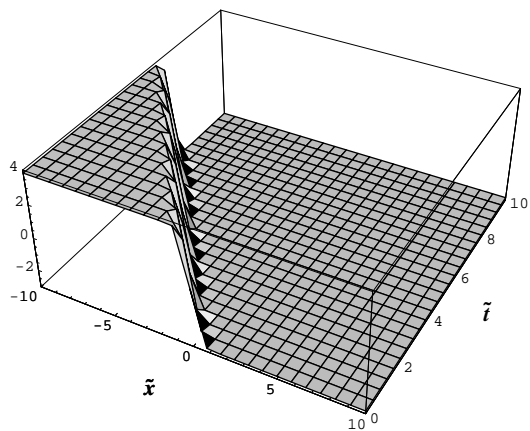}
\caption{\small{The traveling solitary wave for $Re_h=1.9$, $C_1=3$,
$\sigma=\sigma_3 \simeq 0.034 $ }.} \label{Soliton3D}
\end{figure}

\noindent The variations of (\ref{Soliton}) a function of the
non-dimensional space variable $\tilde{x}$ and the non-dimensional
time variable $\tilde{t}$ for $Re_h=1.9$, $C_1=3$,
$\sigma=\sigma_3\simeq 0.034$, is displayed in Figure
\ref{Soliton3D}.
\\

\noindent Analytical calculation yield:

\begin{equation}
\scriptsize{ \tilde {u}_{\tilde {t}}+ c\, \tilde {u}\, \tilde
{u}_{\tilde {x}}-\frac{1}{Re_h}
 \, {\tilde {u}}_{\tilde {x}\tilde
{x}}\tilde{u}(\tilde{x}, \tilde{t}) =\frac{2 \,C_1^2 \,(8 \,C_1^3-9
\,C_1)
   \text{sech}^2\,\big [\,C_1 \,(\tilde {x}-v\,\tilde {t}) \big ]}{\,(8
   \,C_1^2-13)^2 \,Re_h^2 }\,\big [ 8 \,C_1^2-8-18 \,Re_h+2 \,(8
\,C_1^3-9 \,C_1) \,\tanh
   \big [C_1 \,(\tilde {x}-v\,\tilde {t}) \big ]}
\end{equation}

\noindent As $\tilde {t}$ increases, $\scriptsize{ \sup_{\tilde
{t},\,\tilde {x}} \, \mid \tilde {u}_{\tilde {t}}+ c\, \tilde {u}\,
\tilde {u}_{\tilde {x}}-\frac{1}{Re_h}
 \, {\tilde {u}}_{\tilde {x}\tilde
{x}}\tilde{u}(\tilde{x}, \tilde{t})\mid} $ decreases and tends
towards 0; hence, the soliton (\ref{Soliton}) tends towards a
solitary wave solution of (\ref{Burgers}), with an error $o\big (
\frac {e^{-2\,(\tilde {x}-v\,\tilde {t})}}{Re_h^2} \big  )$.

\subsection{Finite difference calculation}

 \noindent Advect $n$ times the soliton (\ref{Soliton}) through the Finite
second-order centered scheme in space, Euler-time scheme
\ref{Scheme1}, in the space domain $x\in [x_0,+\infty[$, for:

\begin{enumerate}
\item [\emph{i}.] $Re_h=1.9$,
$C_1=3$, $\sigma=\sigma_3\simeq 0.034$, which corresponds to the
existence of a solitary wave solution of (\ref{Approx});

\item [\emph{ii}.] $Re_h=1.9$, $\sigma=0.06=\sigma^{'} \neq \sigma_3$, and $Re_h=1.9$, $\sigma=0.07=\sigma^{''}\neq \sigma_3$, which do not correspond to the
existence of a solitary wave solution of (\ref{Approx}).

\end{enumerate}

 \noindent The mesh points number is equal to 150, with 100 points in the front
 wave. $x_0$ is equal to $-20\,\widetilde{h}$.

 \noindent Denote:

\begin{equation}
\widetilde{\tau}={\sigma_3}\,\,\,,
\,\,\,\widetilde{\tau}^{'}={\sigma^{'}}\,\,\,,
\,\,\,\widetilde{\tau}^{''}={\sigma^{''}}
\end{equation}

\begin{figure}[h!]
\centerline{\psfig{file=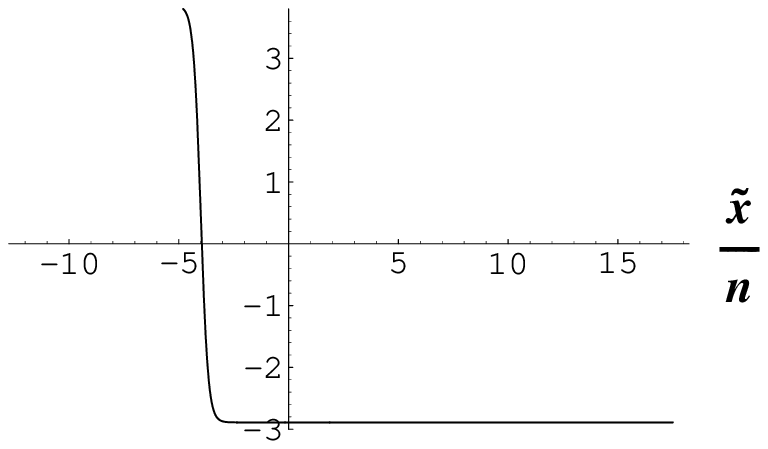,width=7cm}}
\centerline{\psfig{file=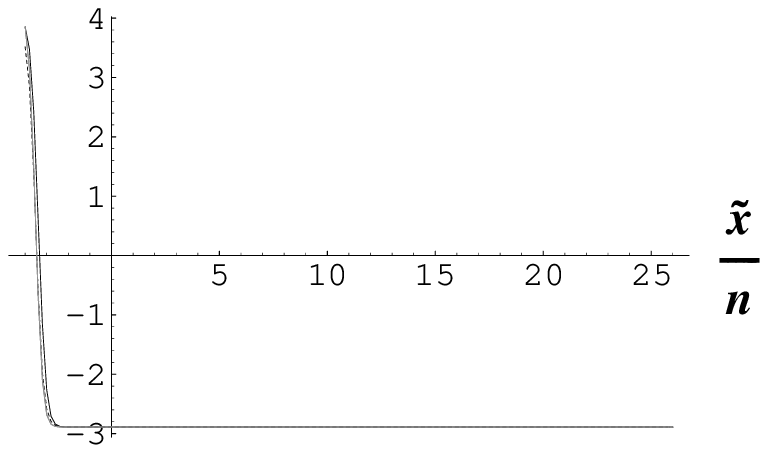,width=7cm}\psfig{file=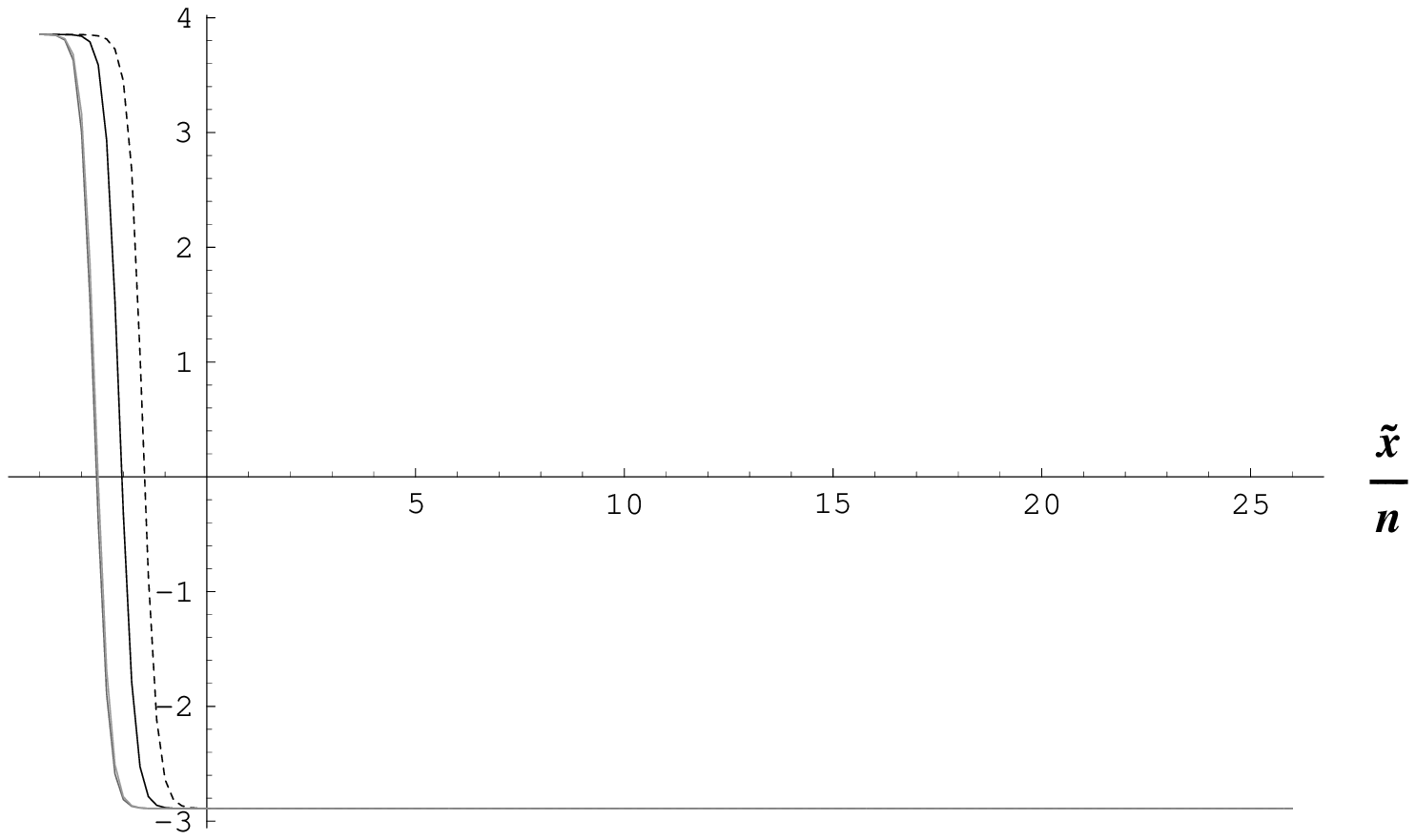,width=7cm}}
\centerline{\psfig{file=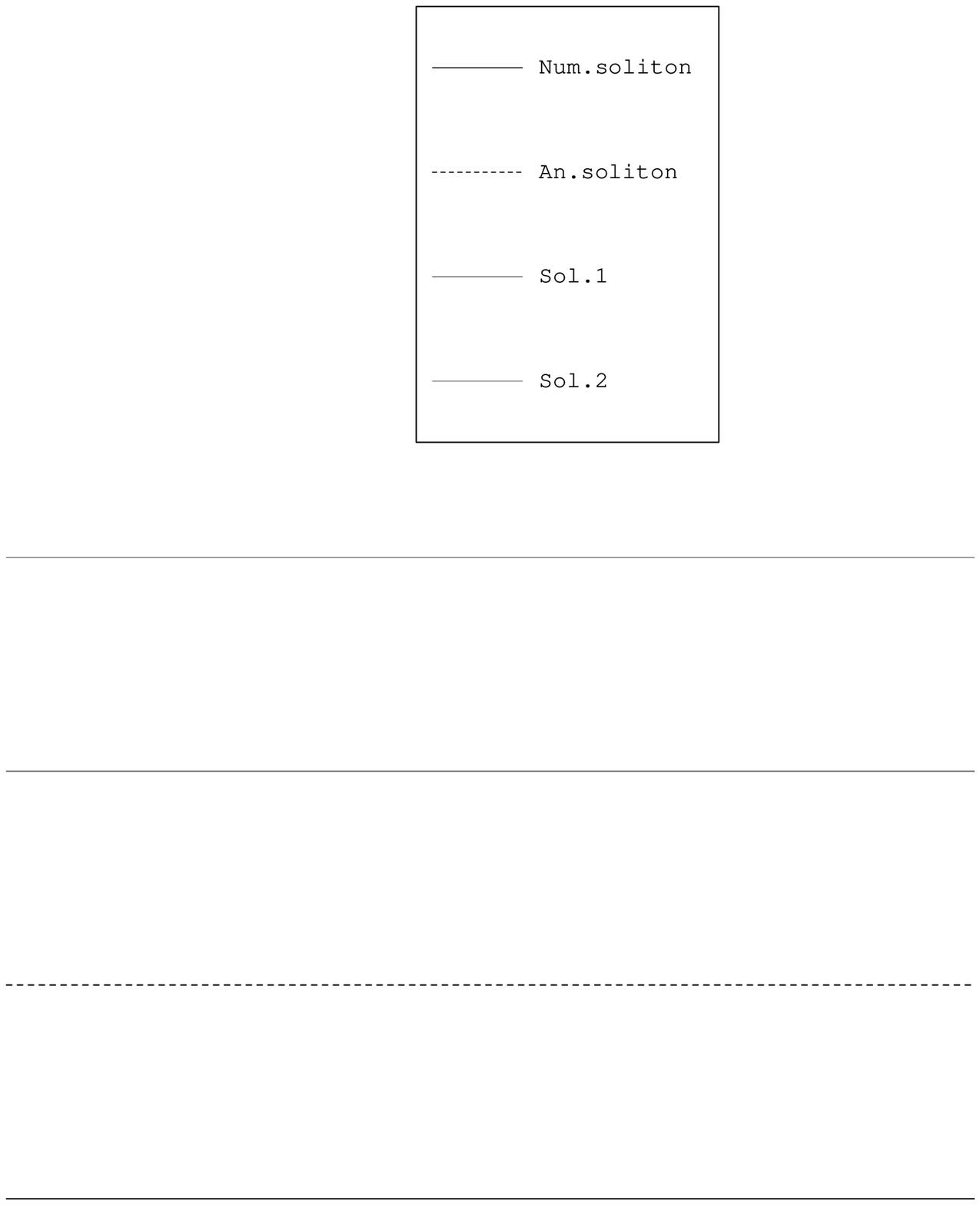,width=3cm}} \caption{Numerical
soliton, analytical soliton, and numerical solutions. Top:
$\widetilde{t}=0$. Bottom: $\widetilde{t}=10\,\widetilde{\tau}$,
$\widetilde{t}=50\,\widetilde{\tau}$.} \label{Solutions}
\end{figure}

 The numerical solitary wave, the analytical solitary wave (\ref{Soliton}), and the numerical solutions
  at $\widetilde {t}=n_t\,\widetilde{\tau}=
  E\big [ \frac {n_t\,\widetilde{\tau}} {\widetilde{\tau}^{'}}\big ]\,\widetilde{\tau}^{'}=
  E\big [ \frac {n_t\,\widetilde{\tau}} {\widetilde{\tau}^{''}}\big ]\,\widetilde{\tau}^{''}$,
  where $E$ denotes the entire part, for $n_t=0$, $n_t=5$, $n_t=10$, $n_t=50$,
  respectively, are displayed in Fig.
\ref{Solutions}.

\newpage

\addcontentsline{toc}{section}{\numberline{}References}

\end{document}